\documentclass[12pt]{article}

\usepackage{amsmath}
\usepackage{amssymb}
\usepackage{amsthm}
\usepackage{mathrsfs, color}
\usepackage{subfigure}
\usepackage{graphicx}
\usepackage{graphics}
\usepackage{pstricks}
\usepackage{pst-grad} 
\usepackage{pst-plot} 
\usepackage{multirow}
\usepackage{lscape}

\newcommand{\ve}[0]{\mathbf e}
\newcommand{\vg}[0]{\mathbf g}

\newcommand{\vs}[0]{\mathbf s}
\newcommand{\vu}[0]{\mathbf u}

\newcommand{\vw}[0]{\mathbf w}
\newcommand{\vx}[0]{\mathbf x}
\newcommand{\vy}[0]{\mathbf y}

\newcommand{\vzero}[0]{\mathbf 0}
\newcommand{\vone}[0]{\mathbf 1}

\newcommand{\Cov}[0]{\text{Cov}}

\newcommand{\diag}[0]{\text{diag}}
\newcommand{\tr}[0]{\text{tr}}

\newcommand{\calG}[0]{\mathcal{G}}
\newcommand{\E}[0]{\mathbb{E}}

\newcommand{\Prob}[0]{\mathbb{P}}
\newcommand{\I}[0]{\mathbb{I}}

\newcommand{\Id}[0]{\text{Id}}

\theoremstyle{plain}
\newtheorem{thm}{Theorem}[section]
\newtheorem{lem}[thm]{Lemma}

\newtheorem{cor}[thm]{Corollary}

\theoremstyle{definition}

\newtheorem{exmp}{Example}[section]

\theoremstyle{remark}
\newtheorem{rmk}{Remark}

\begin{document}

\begin{center}
\textbf{\Large A Note on Moment Inequality for Quadratic Forms} \\ \vspace{0.2in}
\text{Xiaohui Chen}\footnote{Department of Statistics, University of Illinois at Urbana-Champaign. Email: xhchen@illinois.edu. Version: May 7, 2014.} \\
\end{center}

\begin{abstract}
Moment inequality for quadratic forms of random vectors is of particular interest in covariance matrix testing and estimation problems. In this paper, we prove a Rosenthal-type inequality, which exhibits new features and certain improvement beyond the unstructured Rosenthal inequality of quadratic forms when dimension of the vectors increases without bound. Applications to test the block diagonal structures and detect the sparsity in the high-dimensional covariance matrix are presented.
\end{abstract}


\section{Introduction}
\label{sec:introduction}

Covariance matrix plays a central role in multivariate analysis, spatial statistics, pattern recognition and array signal processing. Let $\vx = (X_1,\cdots,X_p)^\top$ be a $p$-variate random vector with mean zero and covariance matrix $\Sigma = \E(\vx \vx^\top)$. Let $\vx_i = (X_{1i},\cdots,X_{pi})^\top$ be independent and identically distributed (iid) copies of $\vx$. To test and estimate the covariance matrix $\Sigma$ based on the observed values of $\vx_i$, there are $p(p+1)/2$ parameters. When the dimension $p$ is large relatively to the sample size $n$, the parameter space grows quadratically in $p$, thus making the covariance testing and estimation problems challenging. Leveraging low-dimensional structures in $\Sigma$, recent literature focuses on structures detection and regularized estimation of the covariance matrix; see e.g. \cite{chenzhangzhong2010a,bickellevina2008a,bickellevina2008b,chenxuwu2013a} among many others. A key step of the success of those regularized estimates relies on sharp large deviation and moment inequalities of the second-order partial sum process involving the quadratic form $\sum_{i=1}^n (X_{ji} X_{ki} - \sigma_{jk})$. Large deviation inequalities can be found in \cite{nagaev1979a,rudelsonvershynin2013a}. In this paper, we shall focus on the moment bounds on the quadratic form, which are useful in controlling the behavior of the second-order statistics of high-dimensional random vectors with heavier tails than the sub-Gaussian distributions.

Classical moment bounds for $\sum_{j=1}^p X_j$ include the well-known Rosenthal's inequality \cite{rosenthal1970}, when $\vx$ has iid components. For the quadratic form, assuming $\E(X_j^2) = 1$, \cite{baisilverstein2009a} established the following inequality
\begin{equation}
\label{eqn:bai-silverstein-quad}
\E|\vx^\top A \vx - \tr(A)|^q \le C_q \left\{ [\kappa_4^4 \tr(A A^\top)]^{q/2} + \kappa_{2q}^{2q} \tr[(A A^\top)^{q/2}] \right\},
\end{equation}
where $A$ is an arbitrary $p \times p$ matrix and $\kappa_w = \|X_j\|_w = (\E|X_j|^w)^{1/w}$ . $C_q$, throughout the paper, denotes a generic constant that only depends on $q$. In essence, (\ref{eqn:bai-silverstein-quad}) is a Rosenthal-type inequality for the product partial sums, weighted by the square matrix $A$. However, we shall see that, in many interesting cases (e.g. Example \ref{exmp:3} and \ref{exam:4} and the applications in Section \ref{sec:application}), (\ref{eqn:bai-silverstein-quad}) is not tight. Therefore, we present a sharper version of moment inequality for bounding $\E|\vx^\top A \vx - \tr(A)|^q$; see Theorem \ref{thm:Rosenthal-quad}. For $w > 0$, a random variable $X \in {\cal L}^w$ iff $\|X\|_w < \infty$. Let $\mathbb{S}^{p \times p}$ be the set of $p \times p$ symmetric matrices and $\mathbb{S}^{p \times p}_+$ be the subset of nonnegative-definite matrices in $\mathbb{S}^{p \times p}$.

\section{Main Result}
\label{sec:main-result}

\begin{thm}
\label{thm:Rosenthal-quad}
Let $X_1, \cdots, X_p$ be iid mean-zero random variables in ${\cal L}^{2q}, q >2$ and $A = (a_{jk})_{1 \le j,k \le p} \in \mathbb{S}^{p \times p}$. Let $\vx = (X_1,\cdots,X_p)^\top$, $s_j^2$ be the eigenvalues of $A^\top A$ and $\vs = (s_1, \cdots, s_p)^\top$. Then, there exists a constant $C_q \in (0, \infty)$ such that
\begin{eqnarray}
\nonumber
\E|\vx^\top A \vx - \tr(A)|^q &\le& C_q \bigg\{ \kappa_{2q}^{2q} (\sum_{j=1}^p a_{jj}^2)^{q/2} + \kappa_{2q}^{2q} (\sum_{1 \le j \neq k \le p} a_{jk}^{2q})^{1/2} \\
\label{eqn:Rosenthal-quad}
& & \qquad + (\kappa_2 \kappa_{2q})^q \left[\sum_{k=1}^p (\sum_{j \neq k} a_{jk}^2)^q\right]^{1/2} + \kappa_2^q |\vs|_4^q \bigg\}.
\end{eqnarray}
\end{thm}

\begin{rmk}
The symmetry of $A$ is not essential in Theorem \ref{thm:Rosenthal-quad}; it is only used for a compact form (\ref{eqn:Rosenthal-quad}). From the proof in Section \ref{sec:proofs}, similar form of (\ref{eqn:Rosenthal-quad}) can be derived for arbitrary square matrix $A$. In addition, symmetry and nonnegative-definiteness of $A$ is a natural requirement for many applications when the quadratic form is the covariance matrix of $\vy = A^{1/2} \vx$. \qed
\end{rmk}

\begin{rmk}
The infinite-dimensional version of (\ref{eqn:Rosenthal-quad}) (and also (\ref{eqn:bai-silverstein-quad})) is also true since the right-hand side of (\ref{eqn:Rosenthal-quad}) is uniform in $p$, provided that the right-hand side converges. This extension allows us to deal with moment bounds of centered covariance process of a linear process of the form $X_j = \sum_{m=0}^\infty b_m \xi_{j-m}$, where $\xi_m$ are iid random variable with mean zero, unit variance and $\xi_i \in {\cal L}^{2q}, q > 2$. As a sufficient condition, we need $b_m = m^{-\beta} \ell(m)$, where $\beta > 1/2$ and $\ell(\cdot)$ is a slowly varying function at $\infty$.\footnote{A slowly varying function $\ell(\cdot)$ in Karamata's sense at $\infty$ is a measurable function such that for any $t > 0$, $\lim_{x\to\infty} \ell(tx) / \ell(x) = 1$. Examples include: constant functions, $\log{x}$, $\exp[(\log{x})^\alpha], 0 < \alpha < 1$.} Note that $\beta \in (1/2, 1)$ allows $X_j$ to have long-range dependence, a widely observed phenomenon in hydrology, network traffic and financial data \cite{giraitiskoulsurgailis2012a}. \qed
\end{rmk}

Now, we compare Theorem \ref{thm:Rosenthal-quad} with the standard moment bound of the quadratic form \cite[Lemma B.26]{baisilverstein2009a}. Since $\tr(A A^\top) = |A|_F^2$ and $\tr[(A A^\top)^{q/2}] = \sum_{j=1}^p |s_j|^q$, (\ref{eqn:bai-silverstein-quad}) is equivalent to
\begin{equation}
\label{eqn:bai-silverstein-quad-alternative}
\E|\vx^\top A \vx - \tr(A)|^q \le C_q (\kappa_4^{2q} |A|_F^q + \kappa_{2q}^{2q} |\vs|_q^q).
\end{equation}
It is interesting to note that both inequalities (\ref{eqn:Rosenthal-quad}) and (\ref{eqn:bai-silverstein-quad-alternative}) balance between the $\ell^2, \ell^4$ and $\ell^q$ norms of the eigenvalues of $A^\top A$, yet in different ways. Let $q$ be fixed, $\kappa_q$ be a constant and the number of variables $p \to \infty$. Then, the Frobenius norm component of our inequality (\ref{eqn:Rosenthal-quad}) is at most the order of (\ref{eqn:bai-silverstein-quad-alternative}) since $(\sum_j a_{jj}^2)^{q/2} + (\sum_{j \neq k} a_{jk}^{2q})^{1/2} + [\sum_k (\sum_{j \neq k} a_{jk}^2)^q]^{1/2} \le 2 |A|_F^q$. On the other hand, since $|\vs|_4 \le |\vs|_q, 2 < q \le 4$ and $|\vs|_4 \ge |\vs|_q, q >4$, the spectral component of (\ref{eqn:Rosenthal-quad}) could be either on smaller or larger order than (\ref{eqn:bai-silverstein-quad-alternative}). In general, it is not possible to have one dominates the other. However, we shall see that, in certain structured cases,  (\ref{eqn:Rosenthal-quad}) improves over (\ref{eqn:bai-silverstein-quad-alternative}) in that the ratio of the two bounds vanishes to zero. Below, we give a few examples.

\begin{exmp}
Consider $A = \Id_{p \times p}$. Then, moment bounds in (\ref{eqn:Rosenthal-quad}) and (\ref{eqn:bai-silverstein-quad}) are $O(p^{q/2})$, both agreeing with the central limit theorem since $\E|\sum_{j=1}^p (X_j^2 - 1)|^q = O(p^{q/2})$.
\end{exmp}

\begin{exmp}
Consider $A = {\mathbf 1}_{p \times p}$, where ${\mathbf 1}_{p \times p}$ denotes the $p \times p$ matrix with all entries equal to one. In this case, moment bounds in (\ref{eqn:Rosenthal-quad}) and (\ref{eqn:bai-silverstein-quad}) both become $O(p^q)$.
\end{exmp}

\begin{exmp}
\label{exmp:3}
Suppose that $m$ and $k$ are two positive integers such that $p = mk, k = o(p)$ and $k, m \to \infty$. Consider the following block diagonal matrix
\begin{equation}
\label{eqn:block-diagonal-example}
A = \left( \begin{array}{cccc}
{\mathbf 1}_{k \times k} & & & \\
& {\mathbf 1}_{k \times k} & & \\
& & {\mathbf 1}_{k \times k} & \\
& & & {\mathbf 1}_{k \times k} \\
\end{array} \right)_{p \times p}.
\end{equation}
This matrix corresponds to a covariance matrix for $p$ variables such that there are $m$ uncorrelated clusters, each of which has $k$ variables perfectly correlated. Then, for this $A$, we have $|A|_F^2 = m k^2, s_1 = \cdots = s_m = k$ and $s_{m+1} = \cdots = s_p = 0$. Thus, $|\vs|_q = m^{1/q} k$ and $|\vs|_4 = m^{1/4} k$. Therefore, the bound (\ref{eqn:bai-silverstein-quad}) is $O(m^{q/2} k^q)$ (also $O(p^{q/2} k^{q/2})$) and our bound (\ref{eqn:Rosenthal-quad}) is $O(p^{q/2} + p^{1/2} k^{q/2} + m^{q/4} k^q)$. Since $m, k \to \infty$, it is clear that $p^{q/2} + p^{1/2} k^{q/2} + m^{q/4} k^q = o(p^{q/2} k^{q/2})$. Therefore, (\ref{eqn:Rosenthal-quad}) is sharper than (\ref{eqn:bai-silverstein-quad}).
\end{exmp}

\begin{exmp}
\label{exam:4}
Let $r \in [0,1)$ and $\rho(A)$ be the spectral norm of a symmetric matrix $A$. Consider the following matrix class
\begin{equation}
\label{eqn:sparse-matrix}
\calG_r(M_p) = \{ A \in \mathbb{S}^{p \times p} : \rho(A) \le C_0, \; \max_{k \le p} \sum_{j=1}^p |a_{jk}|^r \le M_p \},
\end{equation}
where $C_0 > 0$ is a finite constant and $M_p$ is a constant that may depend on $p$. The class $\calG_r(M_p)$ contains the approximately sparse matrices that are invariant under permutation. Here, the sparsity is measured by the size of the strong $\ell^r$-ball. Similar classes of sparse matrices have been extensively studied in the high-dimensional sparse covariance matrix estimation literature; see e.g. \cite{bickellevina2008a,chenxuwu2013a}.

\begin{cor}
\label{cor:sparse-matrix-bound}
Let $q > 2$ and $M_p \to \infty$ as $p \to \infty$. Then, under the conditions in Theorem \ref{thm:Rosenthal-quad}, we have that
\begin{equation}
\label{eqn:sparse-matrix-bound}
\sup_{A \in \calG_r(M_p)} \E|\vx^\top A \vx - \tr(A)|^q \le C(q, r, C_0) (p^{q/2} + p^{1/2} M_p^{q/2}).
\end{equation}
\end{cor}
\end{exmp}

Equation (\ref{eqn:sparse-matrix-bound}) is an application of Theorem \ref{thm:Rosenthal-quad}. On the other hand, for $A \in \calG_r(M_p)$, it is easy to see that $|A|_F^2 \le C_0^{2-r} p M_p$ and $|\vs|_q^q \le C_0^q p$. Therefore, (\ref{eqn:bai-silverstein-quad-alternative}) leads to
\begin{equation*}
\E|\vx^\top A \vx - \tr(A)|^q \le C(q, r, C_0) p^{q/2} M_p^{q/2}.
\end{equation*}
Since $p^{q/2} + p^{1/2} M_p^{q/2} = o(p^{q/2} M_p^{q/2})$ under the condition $M_p \to \infty$, it is clear that the inequality (\ref{eqn:sparse-matrix-bound}) provides a tighter moment bound for the quadratic form with sparsity structure in the matrix $A$.

\section{Applications}
\label{sec:application}

In this section, we provide two applications for testing low-dimensional structures for high-dimensional covariance matrix and its inverse matrix (a.k.a. the precision matrix) for non-Gaussian observations. Let $\vy_1, \cdots, \vy_n$ be $n$ iid mean-zero random vectors such that $\vy_i = (Y_{i1}, \cdots, Y_{ip})^\top$ has independent components with $\E(Y_{ij}^2) = 1$ for $j = 1, \cdots, p$. Let $\Omega \in \mathbb{S}^{p \times p}_+$ and observations $\vx_i = \Omega^{-1/2} \vy_i$ be a linear transformation of $\vy_i$. Then, $\E(\vx_i) = \vzero$ and $\Cov(\vx_i) = \Omega^{-1}$; i.e. the precision matrix of $\vx_i$ is $\Omega$.

\subsection{Test block diagonal covariance and precision matrix}
\label{subsec:app-test-block-diagonal}

Let $A = \diag(A_1, \cdots, A_{m/2})$ and $B = \diag(B_1, \cdots, B_m)$, where $A_i, i = 1,\cdots,m/2,$ are $(2k) \times (2k)$ matrix and $B_i, i = 1,\cdots, m,$ are $k \times k$ matrix. We wish to test the following simple hypotheses for the precision (or equivalently covariance) matrix
\begin{equation}
\label{eqn:precision-matrix-test}
H_0: \Omega = A \qquad \text{V.S.} \qquad H_1: \Omega = B.
\end{equation}
Since the precision matrix fully encodes the partial correlation graph \cite{lauritzen1996a}, (\ref{eqn:precision-matrix-test}) naturally rises in the problems of testing the graph structures of a large number of variables; for instance, different block structures on the pathway and modularity in the genetic regulatory networks; c.f. \cite{tanwittenshojaie2013a,hussamitibshirani2013a}. When $\vx \sim N(0, \Omega^{-1})$, we use the (log-)likelihood ratio test (LRT) statistic
\begin{equation}
\label{eqn:LRT}
L_n = \log\det(B A^{-1}) + {1 \over n} \sum_{i=1}^n \vx_i^\top (A-B) \vx_i.
\end{equation}
If $p$ is fixed and $A-B \in \mathbb{S}^{p \times p}$, then we can show that, under $H_0$, $L^\ast_n = n (L_n - \log\det(B A^{-1}))$ follows a distribution $F$ of the linear combination of $p$ independent $\chi^2(n)$ random variables with coefficients given by the eigenvalues of the matrix $G = A^{-1/2} (A-B) A^{-1/2}$ \footnote{To see this, let $G = P D P^\top$ be its eigen-decomposition and write $L^\ast_n = \sum_{i=1}^n {\vy_i^\top A^{-1/2} (A-B) A^{-1/2} \vy_i} = \sum_{i=1}^n \vw_i^\top D \vw_i$, where $\vw_i = P^\top \vy_i \sim N(0, \Id)$ are iid. Then, $L^\ast_n = \sum_{j=1}^p d_j \chi_j^2$, where $d_j$ are eigenvalues of $G$ and $\chi_j^2$ are iid $\chi^2(n)$ random variables.}. Then, a critical region at the significance level $\alpha$ for the LRT (\ref{eqn:LRT}) is $(r_0, \infty)$ such that $\Prob_F(L^\ast_n > r_0 \mid H_0) = \alpha$. The percentiles of $L^\ast_n$ can be found in \cite{moschopouloscanada1984}. For non-Gaussian $\vx \in {\cal L}^{2q}$, (\ref{eqn:LRT}) can be viewed as a quasi-LRT statistic and construction of the exact critical region is not available, in particular for high-dimensional $\vx_i$. However, a conservative critical region can be constructed as $(\tilde{r}_0, \infty)$ such that
\begin{equation}
\label{eqn:conservative-cr}
\Prob(L_n^\ast - n \tr(G) \ge \tilde{r}_0 \mid H_0) \le \alpha.
\end{equation}
Denote $U_p$ to be the upper bound in (\ref{eqn:Rosenthal-quad}). Since $A$ and $B$ are block diagonal matrices with block size of $(2k) \times (2k)$ and $k \times k$, respectively, it is clear that $G$ is also a block diagonal matrix with block size of $(2k) \times (2k)$. Therefore, by Theorem \ref{thm:Rosenthal-quad} and Markov's inequality, $\tilde{r}_0$ is chosen such that
\begin{equation*}
\Prob ( | \sum_{i=1}^n (\vy_i^\top G \vy_i  - \tr(G)) | \ge \tilde{r}_0 ) \le C_q {n^{q/2} U_p \over \tilde{r}_0^q} \le \alpha.
\end{equation*}
This gives a conservative LRT in the sense of (\ref{eqn:conservative-cr}) and the critical region for $L_n^\ast$ at the significance level $\alpha$ is $(n \tr(G) + n^{1/2} (C_q U_p / \alpha)^{1/q}, \infty)$. As a special case, if $A_i = \vone_{(2k) \times (2k)} + \Id_{(2k) \times (2k)}, i = 1,\cdots,m/2$, and $B = \Id_{p \times p}$, then $G$ has the form (\ref{eqn:block-diagonal-example}). Therefore, the critical region constructed by Theorem \ref{thm:Rosenthal-quad} is more accurate than (\ref{eqn:bai-silverstein-quad}).

\subsection{Sparsity detection in the precision matrix}
\label{subsec:app-test-sparsity}

Here, we present a second application for detecting the sparse graph structures specified by the off-diagonal entries in the precision matrix $\Omega$. There has been a large volume of literature for estimating the sparse graphical models \cite{friedmanhastietibshirani2007a,ravikumar2011a,chenxuwu2013a}. In real data analysis, we would like first to know that whether or not there is a structure in the sparse graph and if so, given the data, how confident? To this end, let $B \in \calG_r(M_p) \setminus \Id_{p \times p}$, we would like to perform the test against the sparse alternative
\begin{equation}
\label{eqn:precision-matrix-test-sparsity}
H_0: \Omega = \Id_{p \times p} \qquad \text{V.S.} \qquad H_1: \Omega = B.
\end{equation}
We use the same definitions for the LRT statistic $L_n$ and $L_n^\ast$ as in Section \ref{subsec:app-test-block-diagonal} with $A = \Id$. Then, $G = \Id -B$ and $G$ has the same off-diagonal entries as $B$ with switched signs. Following the argument in Section \ref{subsec:app-test-block-diagonal}, we can show that
\begin{equation*}
\Prob ( \sum_{i=1}^n (\vy_i^\top G \vy_i  - \tr(G)) \ge \tilde{r}_0 ) \le C_q {n^{q/2} U_p \over \tilde{r}_0^q}.
\end{equation*}
Now, applying Corollary \ref{cor:sparse-matrix-bound}, we obtain a conservative critical region at level $\alpha$ as
 \begin{equation*}
(n \tr(G) + n^{1/2} (p^{1/2} + p^{1/(2q)} M_p^{1/2}) (C_q / \alpha)^{1/q}, \infty).
\end{equation*}
In contrast, if we use (\ref{eqn:bai-silverstein-quad-alternative}), the critical region is $(n \tr(G) + (n p M_p)^{1/2} (C_q / \alpha)^{1/q}, \infty)$. Again, the critical region constructed by Corollary \ref{cor:sparse-matrix-bound} is more accurate than (\ref{eqn:bai-silverstein-quad}).

\section{Proofs}
\label{sec:proofs}

The following version of Rosenthal's inequality, due to \cite{pinelisutev1984a}, for the partial sum of iid random variables is used repeatedly in our proofs.

\begin{lem}
\label{lem:Rosenthal-Pinelis+Utev}
\cite{pinelisutev1984a}. Let $X_i$ be iid mean-zero random variables such that $X_i \in {\cal L}^q$ and $S_n = X_1 + \cdots + X_n$. Then, for $q > 2$,
\begin{equation}
\label{eqn:Rosenthal-Pinelis+Utev}
\E |S_n|^q \le C^q \left[ q^q \sum_{i=1}^n \E|X_i|^q + q^{q/2} (\sum_{i=1}^n \E X_i^2)^{q/2} \right],
\end{equation}
where $C$ is a numeric constant.
\end{lem}

\subsection{Proof of Theorem \ref{thm:Rosenthal-quad}}

Note that
\begin{equation*}
\vx^\top A \vx - \tr(A) = \sum_{j=1}^p a_{jj} (X_j^2 - 1) + \sum_{1 \le j \neq k \le p} a_{jk} X_j X_k.
\end{equation*}
By Burkholder's inequality \cite{rio2009a}, the first term involving the diagonal terms $a_{jj}$ can be bounded by
\begin{equation}
\label{eqn:diagonal}
\|\sum_j a_{jj}(X_j^2 - 1)\|_q^2 \le (q-1) \nu_q^2 \sum_j a_{jj}^2,
\end{equation}
where $\nu_q = \|X_j^2 - 1\|_q$. Next, we deal with the cross-product term. Denote $Q_n = \sum_{1 \le j \neq k \le p} a_{jk} X_j X_k \in {\cal L}^{2q}$. Let $\delta_j$ be iid Bernoulli random variables with success probability $1/2, 1 \le j \le p$ and $\Lambda_\delta = \{\ell \in [p] : \delta_\ell = 1\}$, where $[p] = \{1,\cdots, p\}$. Let $Q_\delta = \sum_{j,k} a_{jk} \delta_j (1-\delta_k) X_j X_k$ and $\E_\delta Q_\delta$ be the expectation taken w.r.t. $\delta$. Then, $Q_n = 4 \E_\delta Q_\delta$ and by Jensen's inequality
\begin{equation*}
\E|Q_n|^{2q} = 4^{2q} \E_X|\E_\delta Q_\delta|^{2q} \le 16^q \E_{X,\delta} |Q_\delta|^{2q}.
\end{equation*}
Observe that we can also write $Q_\delta = \sum_{k \in \Lambda_\delta^c} X_k (\sum_{j \in \Lambda_\delta} a_{jk} X_j)$. By Rosenthal's inequality Lemma \ref{lem:Rosenthal-Pinelis+Utev}, we have that
\begin{equation*}
\E_{(X_k)_{k \in \Lambda_\delta^c}} |Q_\delta|^{2q} \le C^{2q} \left\{ (2q\kappa_{2q})^{2q} \sum_{k \in \Lambda_\delta^c} \left|\sum_{j \in \Lambda_\delta} a_{jk} X_j\right|^{2q} + (2q\kappa_2^2)^q \left[ \sum_{k \in \Lambda_\delta^c} \left( \sum_{j \in \Lambda_\delta} a_{jk} X_j \right)^2 \right]^q \right\},
\end{equation*}
where $C$ is the numeric constant in (\ref{eqn:Rosenthal-Pinelis+Utev}). Therefore, it follows that $\E_X |Q_\delta|^{2q} \le C^{2q} \left[(2q\kappa_{2q})^{2q} \text{I} + (2q\kappa_2^2)^q \text{II}\right]$, where
\begin{equation*}
\text{I} = \sum_{k \in \Lambda_\delta^c} \E_X \left|\sum_{j \in \Lambda_\delta} a_{jk} X_j\right|^{2q}, \qquad \text{II} = \E_X \left[ \sum_{k \in \Lambda_\delta^c} \left( \sum_{j \in \Lambda_\delta} a_{jk} X_j \right)^2 \right]^q.
\end{equation*}
By a second application of Rosenthal's inequality,
\begin{eqnarray}
\nonumber
\text{I} &\le& \sum_{k \in \Lambda_\delta^c} C^{2q} \left[ (2q\kappa_{2q})^{2q} \sum_{j \in \Lambda_\delta} a_{jk}^{2q} + (2q\kappa_2^2)^q (\sum_{j \in \Lambda_\delta} a_{jk}^2)^q \right] \\
\label{eqn:I}
&\le& (2Cq\kappa_{2q})^{2q} \sum_{j \neq k} a_{jk}^{2q} + (2C^2q\kappa_2^2)^q \sum_k \left( \sum_{j \neq k} a_{jk}^2 \right)^q.
\end{eqnarray}
Now, we bound the moment of $\text{II}$ using the idea of decoupling and the reduction principle to Gaussian random variables \cite{rudelsonvershynin2013a}. Let $g_1, \dots, g_p$ be iid $N(0,1)$ and $Z = \sum_{k \in \Lambda_\delta^c} g_k (\sum_{j \in \Lambda_\delta} a_{jk} X_j)$. Then,
\begin{equation*}
Z \mid \delta, (X_j)_{j \in \Lambda_\delta} \sim N(0, \sigma_\delta^2), \qquad \sigma_\delta^2 = \sum_{k \in \Lambda_\delta^c} (\sum_{j \in \Lambda_\delta} a_{jk} X_j)^2.
\end{equation*}
Since $\E_g(Z^{2q}) = \sigma_\delta^{2q} 2^q \pi^{-1/2} \Gamma(q+1/2)$ and $\Gamma(q+1/2) > (q-1)!$,
by Stirling's formula, we see that
\begin{equation*}
\E_{X,g}(Z^{2q}) > c_q  \sqrt{2 \over q} \left({2q \over e}\right)^{q} \; \E_X(\sigma_\delta^{2q}) = c_q \sqrt{2 \over q} \left({2q \over e}\right)^{q} \; \E_X \left[ \sum_{k \in \Lambda_\delta^c} \left( \sum_{j \in \Lambda_\delta} a_{jk} X_j \right)^2 \right]^q.
\end{equation*}
Therefore, to bound the moment of $\text{II}$, it suffices to consider $\E_{X,g}(Z^{2q})$. Conditioned on $g_1, \cdots, g_p$ and noting that $Z = \sum_{j \in \Lambda_\delta} X_j (\sum_{k \in \Lambda_\delta^c} a_{jk} g_k)$, we have that
\begin{equation*}
\E_X(Z_{2q}^{2q}) \le C^{2q} \left\{ (2q\kappa_{2q})^{2q} \sum_{j \in \Lambda_\delta} (\sum_{k \in \Lambda_\delta^c} a_{jk} g_k)^{2q} + (2q\kappa_2^2)^q \left[ \sum_{j \in \Lambda_\delta} (\sum_{k \in \Lambda_\delta^c} a_{jk} g_k)^2 \right]^q \right\}.
\end{equation*}
Taking expectation over $g$ on both sides, we get that $\E_{X,g}(Z^{2q}) \le C^{2q} [(2q\kappa_{2q})^{2q} \text{III} + (2q\kappa_2^2)^q \text{IV}]$, where
\begin{equation}
\label{eqn:II=III+IV}
\text{III} = \sum_{j \in \Lambda_\delta} \E_g (\sum_{k \in \Lambda_\delta^c} a_{jk} g_k)^{2q}, \qquad \text{IV} = \E_g \left[ \sum_{j \in \Lambda_\delta} (\sum_{k \in \Lambda_\delta^c} a_{jk} g_k)^2 \right]^q.
\end{equation}
Since $\sum_{k \in \Lambda_\delta^c} a_{jk} g_k \sim N(0, \sum_{k \in \Lambda_\delta^c} a_{jk}^2)$ and $\Gamma(q+1/2) < q^{-1/2} \Gamma(q+1)$,
\begin{equation*}
\E_g (\sum_{k \in \Lambda_\delta^c} a_{jk} g_k)^{2q} = {2^q \Gamma(q+1/2) \over \sqrt{\pi}} (\sum_{k \in \Lambda_\delta^c} a_{jk}^2)^q < \sqrt{2} \left({2q \over e}\right)^{q} \; (\sum_{k \in \Lambda_\delta^c} a_{jk}^2)^q.
\end{equation*}
Thus,
\begin{equation}
\label{eqn:III}
\text{III} <  \sqrt{2} \left({2q \over e}\right)^{q} \; \sum_{j \in \Lambda_\delta} (\sum_{k \in \Lambda_\delta^c} a_{jk}^2)^q.
\end{equation}
Next, we handle $\text{IV}$. Let $A_\delta = P_\delta A (\Id-P_\delta)$, where $P_\delta$ is the $p \times p$ projection matrix such that $P_\delta \vg = (g_j \I(g_j \in \delta))_{1 \le j \le p}$. Then, $|A_\delta \vg|^2 = \sum_{j \in \Lambda_\delta} (\sum_{k \in \Lambda_\delta^c} a_{jk} g_k)^2$ and $\text{IV} = \E_g |A_\delta \vg|^2$. Let $A_\delta^\top A_\delta = O D O^\top$ be the eigen-decomposition of $A_\delta^\top A_\delta$ where $D$ is a diagonal matrix such that $D = \diag(s_{\delta,1}^2, \cdots, s_{\delta,p}^2)$. Then,
\begin{equation*}
|A_\delta \vg|^2 = \vg^\top A_\delta^\top A_\delta \vg = (O^\top \vg)^\top D (O^\top \vg) \stackrel{d}{=} \vg^\top D \vg = \sum_{j=1}^p s_{j,\delta}^2 g_j^2.
\end{equation*}
Therefore, by Rosenthal's inequality again, we get that
\begin{equation*}
\text{IV} \le C^q \left[q^q \tau_{2q}^{2q} \sum_{j=1}^p s_{j,\delta}^{2q} + q^{q/2} \tau_4^{2q} (\sum_{j=1}^p s_{j,\delta}^4)^{q/2} \right],
\end{equation*}
where $\tau_q = \|g_j\|_q$. Since $P_\delta$ is a projection operator, $s_{j,\delta}^2 \le s_j^2$, we have that
\begin{equation}
\label{eqn:IV}
\text{IV} \le (C q \tau_{2q}^2)^q |\vs|_{2q}^{2q} + (C q^{1/2} \tau_4^2)^q |\vs|_4^{2q} \le C_q |\vs|_4^{2q}.
\end{equation}
Now, (\ref{eqn:Rosenthal-quad}) follows by collecting all terms together and observing that (\ref{eqn:diagonal}), (\ref{eqn:I}), (\ref{eqn:II=III+IV}) (\ref{eqn:III}) and (\ref{eqn:IV}) are all uniform in $\delta$.
\qed

\subsection{Proof of Corollary \ref{cor:sparse-matrix-bound}}

Note that for $A \in \calG_r(M_p)$, $\max_{j \le p} |a_{jj}| \le \rho(A) \le C_0$ and $(\sum_j a_{jj}^2)^{q/2} \le C_0^q p^{q/2}$. Let $A = U \Lambda U^\top = \sum_{\ell=1}^p \lambda_\ell \vu_\ell \vu_\ell^\top$ be the eigen-decomposition of $A$ and $\ve_j = (0,\cdots,0,1,0,\cdots0)^\top$ be the canonical Euclidean basis of $\mathbb{R}^p$. Clearly, $U^\top U = U U^\top = \Id_{p \times p}$ and $|\lambda_j| \le \rho(A), \forall j = 1,\cdots,p$. Then, by the Cauchy-Schwarz inequality, 
\begin{eqnarray*}
|a_{jk}| &=& |\ve_j^\top A \ve_k| \le \sum_{\ell=1}^p |\lambda_j| |\ve_j^\top \vu_\ell \vu_\ell^\top \ve_k| = \sum_{\ell=1}^p |\lambda_j| |u_{j \ell} u_{k \ell}| \\
&\le& (\sum_{\ell=1}^p |\lambda_j| u_{j \ell}^2)^{1/2} (\sum_{\ell=1}^p |\lambda_j| u_{k \ell}^2)^{1/2} \le \rho(A) \le C_0.
\end{eqnarray*}
Therefore,
\begin{equation*}
(\sum_{j < k} a_{jk}^{2q})^{1/2} \le (\sum_j \sum_k |a_{jk}|^r |a_{jk}|^{2q-r})^{1/2} \le C_0^{q-r/2} p^{1/2} M_p^{1/2}.
\end{equation*}
In addition,
\begin{equation*}
[\sum_k (\sum_{j \neq k} a_{jk}^2)^q]^{1/2} \le [\sum_k (C_0^{2-r} \max_{k'} \sum_j |a_{j k'}|^r)^q]^{1/2} \le C_0^{q(1-r/2)} p^{1/2} M_p^{q/2}
\end{equation*}
and $|\vs|_4^q \le C_0^{q/4} p^{q/4}$ since $\max_{j \le p} |s_j| \le C_0$. Now, (\ref{eqn:sparse-matrix-bound}) follows from Theorem \ref{thm:Rosenthal-quad}.
\qed




\end{document}